\documentclass[11pt]{amsart}
\usepackage{graphicx, color}
\usepackage{amscd}
\usepackage{amsmath}
\usepackage{amsfonts}
\usepackage{amssymb}
\usepackage{mathrsfs}

\newtheorem{theorem}{Theorem}

\newtheorem{lemma}[theorem]{Lemma}
\newtheorem{definition}[theorem]{Definition}

\newtheorem{remark}[theorem]{Remark}

\numberwithin{equation}{section}

\renewcommand{\leq}{\leqslant}

\renewcommand{\geq}{\geqslant}
\begin{document}
\title[The Keller-Osserman problem for the k-Hessian operator]{The
Keller-Osserman problem for the k-Hessian operator}
\author[D.-P. Covei]{Dragos-Patru Covei\\Department of Applied Mathematics\\
The Bucharest University of Economic Studies \\
Piata Romana, 1st district, postal code: 010374, postal office: 22, Romania}
\address{}
\email{\texttt{coveid@yahoo.com}}
\keywords{Existence results; Keller-Osserman condition; k-Hessian equation.\\
\phantom{aa} 2010 AMS Subject Classification: Primary: 35J62, 35J92, 35J20;
Secondary: 35J62, 35J15, 47J30.\\
\texttt{coveidragos@yahoo.com}}

\begin{abstract}
A delicate problem is to obtain existence of solutions to the boundary
blow-up elliptic equation%
\begin{equation*}
\sigma _{k}^{1/k}\left( \lambda \left( D^{2}u\right) \right) =g\left(
u\right) \text{ in }\Omega \text{, }\underset{x\rightarrow x_{0}}{\lim }%
u\left( x\right) =+\infty \text{ }\forall x_{0}\in \partial \Omega \text{,}
\end{equation*}%
where $\sigma _{k}^{1/k}\left( \lambda \left( D^{2}u\right) \right) $ is the 
$k$-Hessian operator and $\Omega \subset \mathbb{R}^{N}$ is a smooth bounded
domain. Our goal is to provide a necessary and sufficient condition on $g$
to ensure existence of at least one positive blow-up solution. The main
tools for proving existence are the comparison principle and the method of
sub and supersolutions.
\end{abstract}

\maketitle

%\hfill\today\bigskip

%%%%%%%%%%%%%%%%%%%%%%%%%%%%%%%%%%%%%%%%%%%%%%%%%%%%%%%%%%%%%%%%%%%%%%%

%%%%%%%%%%%%%%%%%%%%%%%%%%%%%%%%%%%%%%%%%%%%%%%%%%%%%%%%%%%%%%%%%%%%%%%%

\section{Introduction}

For $k=1,2,...,N$ define the $k$-Hessian operator as follows 
\begin{equation*}
\sigma _{k}\left( \lambda \left( D^{2}u\right) \right) =\underset{1\leq
i_{1}<...<i_{k}\leq N}{\sum }\lambda _{i_{1}}\cdot ...\cdot \lambda _{i_{k}},
\end{equation*}%
as the $k^{th}$ elementary symmetric polynomial of the Hessian matrix of a $%
C^{2}$ (i.e., a twice continuously differentiable) function $u$ defined over
a bounded smooth domain $\Omega $ (see for details the important works of 
\cite{IIII}-\cite{IF}). Here $\lambda \left( D^{2}u\right) =\left( \lambda
_{1},...,\lambda _{N}\right) $ is the vector of eigenvalues of $D^{2}u$. In
other words, $\sigma _{k}\left( \lambda \left( D^{2}u\right) \right) $ it is
the sum of all $k\times k$ principal minors of the Hessian matrix $D^{2}u$
and so it is a second order differential operator, which may also be called
the $k$-trace of $D^{2}u$ denoted also by 
\begin{equation*}
T_{k}\left[ u\right] :=tr_{k}u_{xx}
\end{equation*}%
where $u_{xx}$ is the Hesse matrix. We would like to mention that 
\begin{equation*}
\sigma _{1}\left( \lambda \left( D^{2}u\right) \right) =\overset{N}{\underset%
{i=1}{\sum }}\lambda _{i}=\Delta u
\end{equation*}%
is the well known classical Laplace operator and 
\begin{equation*}
\sigma _{N}\left( \lambda \left( D^{2}u\right) \right) =\overset{N}{\underset%
{i=1}{\Pi }}\lambda _{i}=\det \left( D^{2}u\right)
\end{equation*}%
is the Monge-Amp\`{e}re operator. Then, for $k\geq 2$ we know that the $k$%
-Hessian operator is a fully nonlinear partial differential operator of
divergence form considered in \cite{E,HI,KR}. The main goal of this paper is
to study the existence of solutions of the following fully nonlinear, second
order partial differential equation with boundary blow-up of the form 
\begin{equation}
\left\{ 
\begin{array}{l}
\sigma _{k}^{1/k}\left( \lambda \left( D^{2}u\right) \right) =\sigma
_{k}^{1/k}\left( \lambda \right) =g\left( u\right) \text{ in }\Omega \text{,}
\\ 
\underset{x\rightarrow x_{0}}{\lim }u\left( x\right) =+\infty \text{ }%
\forall \text{ }x_{0}\in \partial \Omega \text{,}%
\end{array}%
\right.  \label{ma1}
\end{equation}%
where $k\in \left\{ 1,2,...,N\right\} $ and $g:\mathbb{R}\rightarrow \left[
0,\infty \right) $ is a convex function which satisfies:

(G1)\quad $\frac{1}{g}$ is convex on the set where $g>0$, $g^{k}\in
C^{2+\alpha }\left( \mathbb{R},\left[ 0,\infty \right) \right) $ with $%
\alpha \in \left( 0,1\right) $, $g$ is monotone non-decreasing, $g(s)=0$ for
all $s\leq 0$ and $g(s)>0$ for all $s>0$;

(G2)\quad there exists$\mathit{\ }\beta >0\mathit{\ \ }$such that 
\begin{equation*}
\int_{\beta }^{\infty }\frac{1}{\sqrt[k+1]{G(t)-G(\beta )}}dt\in \left(
0,\infty \right) \text{ for\ }G(t)=\int_{0}^{t}g^{k}(z)dz.
\end{equation*}

The problem (\ref{ma1}) belongs to the class of fully nonlinear elliptic
equations and it is closely related to a geometric problem (see \cite{V,VI}
or for more applications \cite{BA,BE,L,RE}). Hence, the $k$-Hessian operator
appears naturally and it is not introduced as a straightforward
generalization of the Laplace or Monge-Amp\`{e}re operator.

The study of existence of blow-up solutions for semilinear elliptic systems
of the form (\ref{ma1}) goes back to the pioneering papers by Osserman \cite%
{O} and Keller \cite{K,KE}. In fact, from the results of \cite{O} and \cite%
{K} we know that, for a given positive, continuous and nondecreasing
function $g$ and a bounded domain $\Omega \subset \mathbb{R}^{N}$, the
semilinear elliptic partial differential equation $\Delta u=g\left( u\right) 
$ in $\Omega $, possesses a blow-up solution $u:\Omega \rightarrow \mathbb{R}
$ if and only if the nowadays called Keller-Osserman condition holds, i. e. 
\begin{equation}
\int_{1}^{\infty }\left( \int_{0}^{t}g\left( s\right) ds\right)
^{-1/2}dt<+\infty .  \label{ko}
\end{equation}

In the present work we will limit ourselves to the development of
mathematical theory for the more general problem (\ref{ma1}). In our
direction, but for the special case $k=1$ or $k=N$, there are many papers
dealing with existence, uniqueness and asymptotic behavior issues for
blow-up solutions of (\ref{ma1}). Here we wish to mention the works of Diaz 
\cite{DC4}, Osserman \cite{O}, Matero \cite{JM,MAT}, Pohozaev \cite{PH} and
Keller \cite{K} (see also references therein). However, excepting the case $%
\Omega =\mathbb{R}^{N}$ studied by Bao-Ji \cite{BAO} and Bao-Ji-Li \cite%
{BAOII}, we don't know any results about the existence of solutions for the
general problem (\ref{ma1}), that so naturally appears in geometry referred
as $k$-Yamabe problem.

Therefore, in contrast to numerous results on the case $k=1$ less is known
about the situation $k\in \left\{ 2,...,N\right\} $. But, a possible
starting point to approach this kind of problems could be works such as \cite%
{BAO}, \cite{DC}, \cite{O}, \cite{MAT}, \cite{K} and \cite{S}.

We begin by stating our result on existence of solutions. Let $\varkappa
_{1},...,\varkappa _{N-1}$ be the set of principal curvatures of $\partial
\Omega $ at $x$ and for $k\in \left\{ 1,2,...,N-1\right\} $ we define the $k$
- curvature $\sigma _{k}\left( \varkappa _{1},...,\varkappa _{N-1}\right) $
of $\partial \Omega $ \ by%
\begin{equation*}
\sigma _{k}\left( \varkappa _{1},...,\varkappa _{N-1}\right) =\sum_{1\leq
i_{1}\leq ...\leq i_{k}\leq N-1}\varkappa _{i_{1}}\cdot ...\cdot \varkappa
_{i_{k}}.
\end{equation*}

\begin{definition}
\label{iiv}Let $k\in \left\{ 2,...,N\right\} $. An open, bounded subset $%
\Omega $ of $\mathbb{R}^{N}$ is said to be $\left( k-1\right) $-convex if $%
\sigma _{i}\geq 0$, for every $x\in \partial \Omega $ and every$\ i\in
\left\{ 1,...,k-1\right\} $.
\end{definition}

\begin{remark}
In particular, the $(N-1)$-convexity for domains is equivalent to the usual
convexity.
\end{remark}

\begin{definition}
\label{iv}Let $k\in \left\{ 2,...,N\right\} $. An open, bounded subset $%
\Omega $ of $\mathbb{R}^{N}$ is said to be stricly $\left( k-1\right) $%
-convex if%
\begin{equation*}
\begin{array}{ll}
i) & \text{it is (}k-1\text{)-convex, i.e. }\sigma _{i}\geq 0\text{, for
every }x\in \partial \Omega \text{ and every}\ i\in \left\{
1,...,k-1\right\} \text{,} \\ 
ii) & \sigma _{k-1}>0\text{ for every }x\in \partial \Omega .%
\end{array}%
\end{equation*}
\end{definition}

Next, we give our main result on existence of solutions\textit{:}

\begin{theorem}
\label{1.C4}Let $~g:\mathbb{R}\rightarrow \left[ 0,\infty \right) $\textit{\
be a\ function satisfying (G1)}. Then, $g$\textit{\ satisfy the
Keller-Osserman type condition (G2) if and only if \ the problem (\ref{ma1})
admits at least one positive solution }$u$\textit{\ in any }bounded,
strictly ($k-1$) - convex and convex domain $\Omega $ in $\mathbb{R}^{N}$ $%
(N\geq 2)$ with $\partial \Omega \in C^{4+\alpha }$ ($\alpha \in \left(
0,1\right) $)\textit{.}
\end{theorem}

Since there is a proof of Theorem \ref{1.C4} for equation (\ref{ma1}) in 
\cite{DC}, \cite{JM} for the case $k=1$, it will be taken to be known in
what follows.

The contributions of our paper are:

\begin{description}
\item[1.] We established a necessary and sufficient condition for the
nonlinearity $g$ such that the considered boundary blowup $k$-Hessian
equation has solution. The sufficient part has been already obtained by
other authors but for particular nonlinearities $g$ that are in $C^{\infty }$
(see for example Colesanti-Salani-Francini \cite{PS} and Salani \cite{S}).

\item[2.] Our methodology is new and it can be applied for more general
nonlinearities depending on the regularity results obtained for the $k$%
-Hessian operator. The necessary part is proved by analyzing the radially
symmetric solution of the considered equation.
\end{description}

The reminder of this paper is organized as follows. The Section \ref{prel}
is devoted to the presentation of some basic results which are needed for
the study of the solutions for (\ref{ma1}). Part of the results will be
fully proven, and, for some of them, only the statements will be exposed.
Section \ref{33} contains the proof of the main result.

\section{Preliminaries \label{prel}}

Since our main Theorem \ref{1.C4} refer to the case when $\Omega $\textit{\ }%
is any bounded, strictly $\left( k-1\right) -$convex and convex domain in $%
\mathbb{R}^{N}$ $(N\geq 2)$ with $\partial \Omega \in C^{4+\alpha }$ we
assume this holds throughout the paper. For $k=1,2,...,N$, denote the set of 
$k$-admissible functions in $\Omega $ by%
\begin{equation*}
\Phi ^{k}\left( \Omega \right) =\left\{ u\in C^{2}\left( \Omega \right)
\left\vert \sigma _{i}\left( \lambda \left( D^{2}u\right) \right) >0\text{
for }i=1,2,...,k\right. \right\} .
\end{equation*}

Next, we collect some auxiliary results. Firstly, we extract from \cite[p. 12%
]{IF} (see also \cite[Theorem 3, p. 264]{CNS}) the following result.

\begin{lemma}
\label{ivf} Let $c\in (0,\infty )$ and $\psi \in C^{2+\alpha }\left( 
\overline{\Omega }\right) $ with $\alpha \in \left( 0,1\right) $. The
Dirichlet problem%
\begin{equation}
\left\{ 
\begin{array}{c}
\sigma _{k}^{1/k}\left( \lambda \left( D^{2}u\right) \right) =\psi >0\text{
in }\overline{\Omega }\text{, }k>1\text{,} \\ 
u=c\text{ on }\partial \Omega%
\end{array}%
\right.  \label{iv54}
\end{equation}%
admits a (unique) admissible solution $u\in \Phi ^{k}\left( \Omega \right)
\cap C^{2}\left( \overline{\Omega }\right) $.
\end{lemma}

Now, an argument similar to \cite[(Definition 2.2, Remark 2.3)]{CIN} leads
to the following definition.

\begin{definition}
Let $g:\mathbb{R}\rightarrow \left[ 0,\infty \right) $ \textit{be a\
function satisfying (G1)} and $c\in (0,\infty )$. A strict subsolution of 
\begin{equation}
\left\{ 
\begin{array}{l}
\sigma _{k}^{1/k}\left( \lambda \left( D^{2}u\right) \right) =g\left(
u\right) \text{ in }\Omega , \\ 
u=c\text{ on }\partial \Omega ,%
\end{array}%
\right.   \label{ma2}
\end{equation}%
is a function $\underline{u}$ from $\Phi ^{k}\left( \Omega \right) $ such
that for some $\delta >0$ we have 
\begin{equation}
\left\{ 
\begin{array}{c}
\sigma _{k}^{1/k}\left( \lambda \left( D^{2}\underline{u}\left( x\right)
\right) \right) \geq g\left( \underline{u}\left( x\right) \right) +\delta 
\text{ for all }x\in \overline{\Omega }, \\ 
\underline{u}=c\text{ on }\partial \Omega .%
\end{array}%
\right.   \label{ma3}
\end{equation}%
Similarly, a supersolution of \ (\ref{ma2}) is a function $\overline{u}$
from $C^{2}\left( \Omega \right) $ such that 
\begin{equation}
\left\{ 
\begin{array}{c}
\sigma _{k}^{1/k}\left( \lambda \left( D^{2}\overline{u}\left( x\right)
\right) \right) \leq g\left( \overline{u}\left( x\right) \right) \text{ for
all }x\in \Omega , \\ 
\overline{u}\geq c\text{ on }\partial \Omega .%
\end{array}%
\right.   \label{ma4}
\end{equation}
\end{definition}

The following variant of the comparison principle will be used. The proof of
the result goes as in \cite{DD} (or: \cite[(Proposition 2.43, p. 187)]{DC4},
Jian \cite[(Lemma 2.3, p. 249)]{FZ}, \cite{H}).

\begin{lemma}
\label{comp} Assume that $w_{1}$, $w_{2}\in C^{2}\left( \Omega \right) \cap
C\left( \overline{\Omega }\right) $ are such that $w_{1}\in \Phi ^{k}\left(
\Omega \right) $. If 
\begin{equation*}
w_{1}\leq w_{2}\text{\ on}\ \partial \Omega ,
\end{equation*}
and 
\begin{equation*}
-\sigma _{k}^{1/k}\left( \lambda \left( D^{2}w_{1}\left( x\right) \right)
\right) \leq -\sigma _{k}^{1/k}\left( \lambda \left( D^{2}w_{2}\left(
x\right) \right) \right) \text{ in }\Omega ,
\end{equation*}%
then $w_{1}\leq w_{2}\ \ $in$\ \ \overline{\Omega }$.
\end{lemma}

The following Lemma can be found in the paper of \cite{S} (see also \cite{G}%
, \cite{EL}, \cite{IT}, \cite{TT} and \cite{T}).

\begin{lemma}
\label{min}Let $c\in (0,\infty )$. Assume that $g:\mathbb{R}\rightarrow %
\left[ 0,\infty \right) $\textit{\ satisfies condition (G1)}. Under these
hypotheses, the problem 
\begin{equation}
\sigma _{k}^{1/k}\left( \lambda \left( D^{2}u\right) \right) =g(u)\;\text{in}%
\;\Omega ,\text{ }u_{\left\vert \partial \Omega \right. }=c  \label{ma7}
\end{equation}%
possesses a unique $k$-admissible solution $u\in \Phi ^{k}\left( \Omega
\right) $ provided that \textit{there exists} a strict subsolution $%
\underline{u}\in \Phi ^{k}\left( \Omega \right) $ and for any $x\in \Omega $
we have $\underline{u}\left( x\right) \leq u\left( x\right) $.
\end{lemma}

Here, we point that Lemma \ref{min} holds true since the comparision
principle given in \cite{G} can be relaxed, see \cite{DD,MMM} for details.

\begin{remark}
If $\Omega =B_{R}$ is a ball from $\mathbb{R}^{N}$ ($N\geq 2$) of radius $R>0
$, then the solution $u$ in Lemma \ref{min} is a radial solution\textit{\
(if it is not so, then we can get another solution by rotating }$u$, but we
have proved that $u$ is the unique solution). Here, we have used the fact
that the $k$-Hessian operator is invariant with respect to rotations.
Moreover, the condition $\frac{1}{g}$ is convex in $\left( 0,\infty \right) $
is not necessary.
\end{remark}

The existence of a subsolution $\underline{u}$ from Lemma \ref{min} and a
supersolution $\overline{u}$ is given in the following Lemma.

\begin{lemma}
\label{subsuper} If $g:\mathbb{R}\rightarrow \left[ 0,\infty \right) $%
\textit{\ satisfies (G1) then:}

i)\quad for any positive constant $c$, the convex positive smallest solution 
$\overline{u}\in C^{2}\left( \Omega \right) $ of the problem 
\begin{equation}
\left\{ 
\begin{array}{l}
\Delta \overline{u}=N\left( C_{N}^{k}\right) ^{-1}g\left( \overline{u}\left(
x\right) \right) \text{, }x\in \Omega , \\ 
\overline{u}\left( x\right) =\infty \text{ on }\partial \Omega ,%
\end{array}%
\right. \text{,}  \label{lapbl}
\end{equation}%
given by the result of \cite[Theorem 3.1., p.1464]{PS}, is a supersolution
of the problem (\ref{ma7}) \textit{and} for any bounded solution $u\in \Phi
^{k}\left( \Omega \right) $ of (\ref{ma7}) we have that $u\left( x\right)
\leq \overline{u}\left( x\right) $ at each point $x\in \Omega $;

ii)\quad the solution $\underline{u}\in \Phi ^{k}\left( \Omega \right) $ of
the problem 
\begin{equation*}
\left\{ 
\begin{array}{c}
\sigma _{k}^{1/k}\left( \lambda \left( D^{2}\underline{u}\right) \right)
=g\left( c\right) +1\text{ in }\overline{\Omega }\text{, }k>1\text{,} \\ 
\underline{u}=c\text{ on }\partial \Omega ,%
\end{array}%
\right.
\end{equation*}%
given by Lemma \ref{ivf}, is a strict subsolution of the problem (\ref{ma7})
for any positive constant $c$;

iii)\quad the strict subsolution $\underline{u}\in \Phi ^{k}\left( \Omega
\right) $ and the convex supersolution $\overline{u}\in C^{2}\left( \Omega
\right) $ determined in i) and ii) are such that $\underline{u}\leq 
\overline{u}$ in $\overline{\Omega }$.
\end{lemma}

\subparagraph{Proof}

i) The Maclaurin's inequalities%
\begin{equation}
\frac{1}{N}\Delta \overline{u}\geq \left[ \left( C_{N}^{k}\right)
^{-1}\sigma _{k}\left( \lambda \left( D^{2}\overline{u}\right) \right) %
\right] ^{1/k}\text{ for any }k=1,...,N\text{,}  \label{new}
\end{equation}%
where%
\begin{equation*}
C_{N}^{k}=\frac{(N-1)!}{k!(N-k)!}
\end{equation*}%
is the binomial coefficient, gives%
\begin{eqnarray*}
N\left[ \left( C_{N}^{k}\right) ^{-k}\sigma _{k}\left( \lambda \left( D^{2}%
\overline{u}\right) \right) \right] ^{1/k} &\leq &N\left[ \left(
C_{N}^{k}\right) ^{-1}\sigma _{k}\left( \lambda \left( D^{2}\overline{u}%
\right) \right) \right] ^{1/k} \\
&\leq &\Delta \overline{u}=N\left( C_{N}^{k}\right) ^{-1}g\left( \overline{u}%
\left( x\right) \right) .
\end{eqnarray*}%
Thus%
\begin{equation*}
\sigma _{k}^{1/k}\left( \lambda \left( D^{2}\overline{u}\right) \right) \leq
g\left( \overline{u}\left( x\right) \right) \text{ in }\Omega \text{.}
\end{equation*}%
The rest of the proof is a consequence of Lemma \ref{comp}.

ii) To prove the affirmation we observe that \ 
\begin{equation*}
\left\{ 
\begin{array}{c}
\sigma _{k}^{1/k}\left( \lambda \left( D^{2}\underline{u}\right) \right)
=g\left( c\right) +1\geq g\left( \underline{u}\right) \text{ in }\Omega , \\ 
\underline{u}\left( x\right) _{\left\vert \partial \Omega \right.
}=c_{\left\vert \partial \Omega \right. }=c,%
\end{array}%
\right.
\end{equation*}%
as the author \cite{S} observed.

iii) Again, we use the Maclaurin's inequalities 
\begin{equation*}
\frac{1}{N}\Delta \underline{u}\geq \left[ \left( C_{N}^{k}\right)
^{-1}\sigma _{k}\left( \lambda \left( D^{2}\underline{u}\right) \right) %
\right] ^{1/k}\text{ in }\Omega \text{, for any }k=1,...,N\text{,}
\end{equation*}%
to see that%
\begin{equation}
\frac{1}{N}\Delta \underline{u}\geq \left[ \left( C_{N}^{k}\right)
^{-1}\left( g\left( \underline{u}\right) +1\right) ^{k}\right] ^{1/k}\text{
in }\Omega .  \label{ns}
\end{equation}%
Suppose by counterposition that the open set $\omega =\left\{ x\in \Omega
\left\vert \underline{u}>\overline{u}\right. \right\} $ is non empty.
Without loss of generality, we can therefore assume that $\omega $ is
connected, also we work with a connected component of $\omega $. Then (\ref%
{ns}) becomes%
\begin{eqnarray*}
\frac{1}{N}\Delta \underline{u} &\geq &\left[ \left( C_{N}^{k}\right)
^{-1}\left( g\left( \underline{u}\right) +1\right) ^{k}\right] ^{1/k}\geq %
\left[ \left( C_{N}^{k}\right) ^{-1}g^{k}\left( \overline{u}\right) \right]
^{1/k} \\
&\geq &\left[ \left( C_{N}^{k}\right) ^{-k}g^{k}\left( \overline{u}\right) %
\right] ^{1/k}=\frac{1}{N}\Delta \overline{u}\text{ in }\omega ,
\end{eqnarray*}%
from which we have $\Delta \left( \underline{u}-\overline{u}\right) \geq 0$
in $\omega $ and therefore, by the classical maximum principle for the
laplacian, we obtain that $\underline{u}-\overline{u}\leq 0$ in $\omega $.
This is a contradiction with the assumption. The proof is now completed.

For the readers' convenience, we recall the radial form of the $k$-Hessian
operator.

\begin{remark}
(see \cite[(Lemma 2.1, p. 178)]{BAO})For $R>0$, let $B_{R}\left( 0\right)
=\left\{ x\in \Omega \left\vert \left\vert x\right\vert <r\right. \right\} $%
. Assume $\varphi \in C^{2}\left( \left[ 0,R\right) \right) $, with $\varphi
^{\prime }\left( 0\right) =0$. If $u\left( x\right) =\varphi \left( r\right) 
$, where $r=\left\vert x\right\vert <R$, we have that $u\in C^{2}\left(
B_{R}\left( 0\right) \right) $, 
\begin{eqnarray*}
\lambda \left( D^{2}u\left( r\right) \right) &=&\left\{ 
\begin{array}{l}
\left( \varphi ^{\prime \prime }\left( r\right) ,\frac{\varphi ^{\prime
}\left( r\right) }{r},...,\frac{\varphi ^{\prime }\left( r\right) }{r}%
\right) \text{ for }r\in \left( 0,R\right) , \\ 
\left( \varphi ^{\prime \prime }\left( 0\right) ,\varphi ^{\prime \prime
}\left( 0\right) ,...,\varphi ^{\prime \prime }\left( 0\right) \right) \text{
for }r=0%
\end{array}%
\right. \\
\sigma _{k}\left( \lambda \left( D^{2}u\left( r\right) \right) \right)
&=&\left\{ 
\begin{array}{l}
C_{N-1}^{k-1}\varphi ^{\prime \prime }(r)\left( \frac{\varphi ^{\prime }(r)}{%
r}\right) ^{k-1}+C_{N-1}^{k-1}\frac{N-k}{k}\left( \frac{\varphi ^{\prime }(r)%
}{r}\right) ^{k}\text{, }r\in \left( 0,R\right) , \\ 
C_{N}^{k}\left( \varphi ^{\prime \prime }\left( 0\right) \right) ^{k}\text{
for }r=0,%
\end{array}%
\right.
\end{eqnarray*}%
where the prime denotes differentiation with respect to $r=\left\vert
x\right\vert $ and $k\in \left\{ 1,2,...,N\right\} $.
\end{remark}

The following Lemma is a consequence of results from a number of works, we
mention \cite{DC}, \cite[Lemma 2.1, p. 60]{US}.

\begin{lemma}
\label{echiv} Assume that $g:\mathbb{R}\rightarrow \left[ 0,\infty \right) $
satisfies condition (G1). The following statements are equivalent

KO1)\quad there \textit{exists }$\beta >0$\textit{\ \ such that\ }%
\begin{equation}
\text{\ }\mathcal{K}(\beta ):=\int_{\beta }^{\infty }\frac{1}{\sqrt[k+1]{%
G(t)-G(\beta )}}dt\in \left( 0,\infty \right) \text{, }  \label{ma8}
\end{equation}

KOL1)\quad the Sharpened Keller-Osserman condition\quad\ \textit{\ }%
\begin{equation}
\underset{\beta \rightarrow \infty }{\lim }\inf \mathcal{K}(\beta )=0.
\label{ma9}
\end{equation}
\end{lemma}

The next estimates is almost identical to that of \cite{DC}.

\begin{lemma}
\label{2.4C4} Assume that $g:\mathbb{R}\rightarrow \left[ 0,\infty \right) $%
\textit{\ satisfies (G1)}. \textit{If }$\xi \in C^{2}(0,\rho )$\textit{\ is
a non-decreasing function solving}%
\begin{equation}
C_{N-1}^{k-1}\left[ \frac{r^{N-k}}{k}\left( \xi ^{^{\prime }}(r)\right) ^{k}%
\right] ^{\prime }=r^{N-1}g^{k}(\xi \left( r\right) )\text{ in }(0,\rho ),
\label{ma10}
\end{equation}%
then\textit{, for }$0$\textit{\ }$<\rho _{1}<\rho _{2}<\rho ,$ we have that%
\begin{equation}
\int_{\xi (\rho _{1})}^{\xi (\rho _{2})}\frac{\left( C_{N-1}^{k-1}\right)
^{1/\left( k+1\right) }}{\sqrt[k+1]{\left( k+1\right) G\left( r\right)
-G\left( \xi (\rho _{1})\right) }}dr\geq \frac{k}{N-2k}\rho _{1}^{\frac{2k}{%
k+1}}\left[ 1-\left( \frac{\rho _{1}}{\rho _{2}}\right) ^{\frac{N}{k}-2}%
\right] ,  \label{ma11}
\end{equation}%
given that\textit{\ }$N\neq 2k,$\textit{\ and }%
\begin{equation}
\int_{\xi (\rho _{1})}^{\xi (\rho _{2})}\frac{\left( C_{N-1}^{k-1}\right)
^{1/\left( k+1\right) }}{\sqrt[k+1]{\left( k+1\right) G\left( r\right)
-G\left( \xi (\rho _{1})\right) }}dr\geq \rho _{1}^{\frac{2k}{k+1}}\ln \frac{%
\rho _{2}}{\rho _{1}},  \label{ma12}
\end{equation}%
\textit{if }$N=2k$\textit{.}
\end{lemma}

\subparagraph{Proof}

A simple calculation show that for $r\in (\rho _{1},\rho _{2})$ we have that 
$\xi ^{\prime }(r)>0$, see \cite[p. 286]{S}. Moreover, (\ref{ma10}) is
equivalent to 
\begin{equation}
C_{k-1}^{N-1}\xi ^{\prime \prime }(r)\left( \frac{\xi ^{\prime }(r)}{r}%
\right) ^{k-1}+C_{N-1}^{k-1}\frac{N-k}{k}\left( \frac{\xi ^{\prime }(r)}{r}%
\right) ^{k}=g^{k}(\xi (r)).  \label{ma13}
\end{equation}%
Multiplying the equation (\ref{ma13}) by $r^{N+\frac{N}{k}-2}\xi ^{\prime
}(r)$, we get that%
\begin{equation}
\left[ \left( r^{\frac{N}{k}-1}\xi ^{\prime }(r)\right) ^{k+1}\right]
^{\prime }=\frac{k+1}{C_{N-1}^{k-1}}g^{k}(\xi (r))r^{N+\frac{N}{k}-2}\xi
^{\prime }(r).  \label{ma14}
\end{equation}%
Integrating (\ref{ma14}) from $\rho _{1}$ to $r$ we obtain%
\begin{eqnarray}
\left( r^{\frac{N}{k}-1}\xi ^{\prime }(r)\right) ^{k+1} &\geq &\left( r^{%
\frac{N}{k}-1}\xi ^{\prime }(r)\right) ^{k+1}-\left( \rho _{1}^{\frac{N}{k}%
-1}\xi ^{\prime }(\rho _{1})\right) ^{k+1}  \notag \\
&=&\int_{\rho _{1}}^{r}\frac{k+1}{C_{N-1}^{k-1}}g^{k}(\xi (s))s^{N+\frac{N}{k%
}-2}\xi ^{\prime }(s)ds  \label{above} \\
&\geq &\frac{k+1}{C_{N-1}^{k-1}}\rho _{1}^{N+\frac{N}{k}-2}\left( G\left(
\xi (r)\right) -G\left( \xi (\rho _{1})\right) \right) .  \notag
\end{eqnarray}%
Especially, by (\ref{above}) we have that%
\begin{equation*}
r^{\frac{N}{k}-1}\xi ^{\prime }(r)\geq \left( \frac{k+1}{C_{N-1}^{k-1}}%
\right) ^{1/\left( k+1\right) }\rho _{1}^{\frac{N}{k}-\frac{2}{k+1}}\sqrt[k+1%
]{G\left( \xi (r)\right) -G\left( \xi (\rho _{1})\right) }.
\end{equation*}%
Equivalently, this can be written in the following way%
\begin{equation}
\left( \frac{C_{N-1}^{k-1}}{k+1}\right) ^{1/\left( k+1\right) }\frac{\xi
^{^{\prime }}(r)}{\sqrt[k+1]{G\left( \xi (r)\right) -G\left( \xi (\rho
_{1})\right) }}\geq (\frac{\rho _{1}}{r})^{\frac{N}{k}-1}\rho _{1}^{1-\frac{2%
}{k+1}}.  \label{ma15}
\end{equation}%
Integrating the relation (\ref{ma15}) from $\rho _{1}$ and $\rho _{2}$, we
obtain:

\textit{given that }$2k\neq N$ the relation holds%
\begin{eqnarray*}
&&\left( \frac{C_{N-1}^{k-1}}{k+1}\right) ^{1/\left( k+1\right) }\int_{\xi
(\rho _{1})}^{\xi (\rho _{2})}\frac{1}{\sqrt[k+1]{G\left( r\right) -G\left(
\xi (\rho _{1})\right) }}dt \\
&\geq &\int_{\rho _{1}}^{\rho _{2}}(\frac{\rho _{1}}{t})^{\frac{N}{k}-1}\rho
_{1}^{1-\frac{2}{k+1}}dt=\rho _{1}^{\frac{N}{k}-\frac{2}{k+1}}\left( \left. 
\frac{t^{2-\frac{N}{k}}}{2-\frac{N}{k}}\right\vert _{\rho _{1}}^{\rho
_{2}}\right) \\
&=&\frac{k}{N-2k}\rho _{1}^{\frac{N}{k}-\frac{2}{k+1}}\left( \rho _{1}^{2-%
\frac{N}{k}}-\rho _{2}^{2-\frac{N}{k}}\right) ,
\end{eqnarray*}

and the inequality holds%
\begin{eqnarray*}
\left( \frac{C_{N-1}^{k-1}}{k+1}\right) ^{1/\left( k+1\right) }\int_{\xi
(\rho _{1})}^{\xi (\rho _{2})}\frac{1}{\sqrt[k+1]{G\left( r\right) -G\left(
\xi (\rho _{1})\right) }}dr &\geq &\rho _{1}^{1-\frac{2}{k+1}}\int_{\rho
_{1}}^{\rho _{2}}\frac{\rho _{1}}{t}dt \\
&=&\rho _{1}^{2-\frac{2}{k+1}}\ln \frac{\rho _{2}}{\rho _{1}},
\end{eqnarray*}

\textit{if} $2k=N$. Thus we get the conclusion.

\bigskip We can also obtain an estimate as in (\ref{ma11})-(\ref{ma12})
using Maclaurin's inequalities (\ref{new}).

An important consequence of (\ref{ma11}) and (\ref{ma12}) is the next:

\begin{lemma}
\label{bila} Let $g:\mathbb{R}\rightarrow \left[ 0,\infty \right) $\textit{\
be a function satisfying (G1)}. We have: $g$ \textit{satisfy the
Keller-Osserman type condition (G2) if and only if the problem (\ref{ma1})
admits at least one solution }$u\in \Phi ^{k}\left( B_{\rho }\right) $%
\textit{\ on some ball }$B_{\rho }$\textit{.}
\end{lemma}

\subparagraph{Proof}

\bigskip We deal with the first implication. To prove it we shall proceed as
follows. If $2k\neq N$ we assume temporarily that 
\begin{equation*}
\left( C_{N-1}^{k-1}\right) ^{1/\left( k+1\right) }\mathcal{K}(\beta )<\frac{%
k}{\left\vert N-2k\right\vert }.
\end{equation*}%
Let $\underline{\xi }\in \Phi ^{k}\left( B_{1}\right) $ be the strict
subsolution constructed in Lemma \ref{subsuper} with $\Omega =B_{1}$, $%
c=\beta $ and $\overline{\xi }\in C^{2}\left( B_{1}\right) $ be the
supersolution constructed in that way. It is clear that 
\begin{equation*}
\underline{\xi }\left( x\right) \leq \overline{\xi }\left( x\right) \text{
for }x\in \overline{\Omega }.
\end{equation*}%
An alternative strict sub and supersolution in the space $\Phi ^{k}\left(
B_{1}\right) $ can be obtained as in \cite{BAO,BAOII}. From \textit{Lemma %
\ref{min}} in connection with Lemma \ref{subsuper} we know that there exists
a unique radial solution $\xi \in C^{2}\left( B_{1}\right) $ that solves the
problem (\ref{ma2}) with $\Omega =B_{1}$ and such that $\underline{\xi }\leq
\xi \leq \overline{\xi }$. By the knowledge of the classical theory for
ordinary differential equations (see, for example \cite[(Lemma 9, p. 2148)]%
{BAOII} or \cite{BAO}), we know that choosing $\widetilde{\beta }=\xi (0)$
and $\xi ^{\prime }(0)=0,$ the solution $\xi (r):=u(x)$ (for $r=\left\vert
x\right\vert $) can be extended to maximal interval $\left[ 0,\rho \right) $%
. Let us point that in the paper of \cite{BAOII} we have all the discussion
to understand our problem. \ Assume that $\rho <\infty $, then $u$ is an
blow-up solution in the ball $B_{\rho }$. Indeed, by the definition of $\rho 
$, we have 
\begin{equation*}
\text{or }\xi (\rho )=+\infty \text{ or }\xi ^{\prime }(\rho )=+\infty .
\end{equation*}%
In the case $\xi ^{\prime }(\rho )=+\infty $, integrating between $0$ and $r$
in (\ref{ma14})\ we obtain that%
\begin{equation}
\left( r^{\frac{N}{k}-1}\xi ^{\prime }(r)\right) ^{k+1}=\frac{k+1}{%
C_{N-1}^{k-1}}\left[ G(\xi (r))r^{N+\frac{N}{k}-2}-\left( N+\frac{N}{k}%
-2\right) \int_{0}^{r}G\left( \xi (s\right) s^{N+\frac{N}{k}-3}ds\right] .
\label{ma16}
\end{equation}%
Using the equality (\ref{ma16}) we get%
\begin{equation}
r^{\frac{N}{k}-1}\xi ^{\prime }(r)\leq \left( \frac{k+1}{C_{N-1}^{k-1}}%
\right) ^{1/\left( k+1\right) }G^{1/\left( k+1\right) }(\xi (r))r^{\frac{N}{k%
}-\frac{2}{k+1}}\text{, }r\in \left[ 0,\rho \right) .  \label{ma16s}
\end{equation}%
We write Eq. (\ref{ma16s}) in the form%
\begin{equation*}
\xi ^{\prime }(r)\leq \left( \frac{k+1}{C_{N-1}^{k-1}}\right) ^{1/\left(
k+1\right) }G^{1/\left( k+1\right) }(\xi (r))r^{\frac{k-1}{k+1}}\text{, }%
r\in \left[ 0,\rho \right) .
\end{equation*}%
Then $G(\xi (\rho ))=+\infty $, $\xi (\rho )=+\infty $ and as a conclusion $%
\xi (r):=u(x)$ is an blow-up solution. We next turn to the proof of $\rho
<\infty $. In contrary, using \textit{Lemma\ \ref{2.4C4}} with $\ \rho
_{1}=1 $ and $\rho _{2}>1,$ we observe that 
\begin{equation*}
\left( C_{N-1}^{k-1}\right) ^{1/\left( k+1\right) }\mathcal{K}(\beta )\geq 
\frac{k}{N-2k}\left[ 1-\left( \frac{1}{\rho _{2}}\right) ^{\frac{N}{k}-2}%
\right] ,
\end{equation*}%
if $N\neq 2k$ and 
\begin{equation*}
\left( C_{N-1}^{k-1}\right) ^{1/\left( k+1\right) }\mathcal{K}(\beta )\geq
\ln \rho _{2},
\end{equation*}%
if $N=2k$. Using the \textit{Lemma \ref{echiv}} for $\rho _{2}$ approaching $%
\infty $, we obtain a contradiction if either $N=2k$ or 
\begin{equation*}
\left( C_{N-1}^{k-1}\right) ^{1/\left( k+1\right) }\mathcal{K}(\beta )<\frac{%
k}{\left\vert N-2k\right\vert }.
\end{equation*}%
In the case $N\neq 2k$ and $\mathcal{K}(\beta )\geq \frac{k}{\left\vert
N-2k\right\vert }$, direct calculation prove that we can choose $c_{1}>0$
sufficiently large such that%
\begin{equation*}
\frac{1}{c_{1}}\left( C_{N-1}^{k-1}\right) ^{1/\left( k+1\right) }\mathcal{K}%
(\beta )<\frac{k}{\left\vert N-2k\right\vert }.
\end{equation*}%
Repeating the above discussion we obtain that if $g$ is replaced by a well
determined constant $c_{1}^{\left( k+1\right) /k}$ multiplied with $g$, then 
$\widetilde{u}(x):=u(x/c_{1})$ is an blow-up solution of the problem (\ref%
{ma10}) with the nonlinear function $g$ in $B_{\rho c_{1}}$ which denote the
concentric ball of radius $\rho c_{1}$. We have already checked that if the
Keller-Osserman type condition holds,\ then there exists some ball in which
the solution $u$ blows up to finite value of $\rho $.

Our next step is to prove the second implication. For this, we assume that
there exists some ball $B$ of radius $\rho $, whose center we may always
assume to be the origin, in which\ $u(x)$ solves the problem (\ref{ma1}). By
the above proof and \cite{BAOII}, we may always assume that $u$ is a radial
solution and we define $\xi (r)=u(x)$ for $r=\left\vert x\right\vert $, so
that $\xi $ verify the problem (\ref{ma10}) in $\left[ 0,\rho \right) $. A
calculation akin to that in \textit{Lemma \ref{2.4C4}} leads to%
\begin{equation*}
\left[ \left( r^{\frac{N}{k}-1}\xi ^{\prime }(r)\right) ^{k+1}\right]
^{\prime }=\frac{k+1}{C_{N-1}^{k-1}}g^{k}(\xi (r))r^{N+\frac{N}{k}-2}\xi
^{\prime }(r).
\end{equation*}%
Integrating this equation from $0$ to $r$, we have%
\begin{eqnarray*}
\left( r^{\frac{N}{k}-1}\xi ^{\prime }(r)\right) ^{k+1} &=&\int_{0}^{r}\frac{%
k+1}{C_{N-1}^{k-1}}g^{k}(\xi (z))z^{N+\frac{N}{k}-2}\xi ^{\prime }(z)dz \\
&\leq &\frac{k+1}{C_{N-1}^{k-1}}r^{N+\frac{N}{k}-2}[G(\xi (r))-G(\xi (0))].
\end{eqnarray*}%
Evidently, 
\begin{equation*}
r^{\frac{N}{k}-1}\xi ^{\prime }(r)\leq \left( \frac{k+1}{C_{N-1}^{k-1}}%
\right) ^{1/\left( k+1\right) }r^{\frac{N}{k}-\frac{2}{k+1}}\sqrt[k+1]{G(\xi
(r))-G(\xi (0))}.
\end{equation*}%
Integrating once more between $0$ and $\rho $, it follows that%
\begin{equation}
0\leq \int_{0}^{\rho }\frac{\xi ^{^{\prime }}(r)}{\sqrt[k+1]{\left(
k+1\right) \left( G(\xi (r))-G(\xi (0))\right) }}dr\leq \frac{\left(
C_{N-1}^{k-1}\right) ^{-1/\left( k+1\right) }\left( k+1\right) }{2k}\rho ^{%
\frac{2k}{k+1}}.  \label{ma17}
\end{equation}%
So the conclusion can be obtained by choosing $\beta =\xi (0)$ in (\ref{ma8}%
), which finishes the proof of the lemma.

The following result shows the existence of solutions on small balls.

\begin{lemma}
\label{2.7C4}Assume that $g:\mathbb{R}\rightarrow \left[ 0,\infty \right) $
satisfies condition (G1). If (\ref{ma9}) holds and $B_{\rho }$\textit{\ is a
ball of radius }$\rho $\textit{\ then }%
\begin{equation}
\rho _{0}=\inf \{\rho >0:\text{(\ref{ma1}) has a solution in }\Omega
=B_{\rho }\}=0\text{.}  \label{bestC4}
\end{equation}
\end{lemma}

It is an easy exercise to prove Lemma \ref{2.7C4}. To do this we use the
results in \cite{BAOII} for the $k$-Hessian operator.

\subparagraph{Proof}

We assume by contradiction that $\rho _{0}>0$. Let $\beta _{n}$ be a
sequence of real numbers increasing to infinity and satisfying%
\begin{equation*}
\underset{\beta _{n}\rightarrow \infty }{\lim }\text{\ }\inf \mathcal{K}%
\left( \beta _{n}\right) =0.
\end{equation*}%
Let $\underline{u}=\underline{u}^{n}\in \Phi ^{k}\left( B_{\rho
_{0}/2}\right) $ be the radial strict subsolution obtained in Lemma \ref{min}
in connection with Lemma \ref{subsuper} corresponding to $\Omega =B_{\rho
_{0}/2}$, $c=\beta _{n}$ and $\overline{u}$ the supersolution constructed.
It follows in a standard fashion, that for $\underline{u}$ there exists a
unique radial solution in $C^{2}\left( B_{\rho _{0}/2}\right) $, denoted by $%
u_{n}$, of the problem 
\begin{equation*}
\left\{ 
\begin{array}{ll}
\sigma _{k}^{1/k}\left( \lambda \left( D^{2}u_{n}\right) \right) =g(u_{n}) & 
\text{\textit{in} }B_{\rho _{0}/2}, \\ 
u_{n}=\beta _{n} & \text{\textit{on} }\partial B_{\rho _{0}/2}.%
\end{array}%
\right.
\end{equation*}%
Using the definition of $\rho _{0}$ and considering the problem 
\begin{equation}
\left\{ 
\begin{array}{l}
C_{N-1}^{k-1}\left[ \frac{r^{N-k}}{k}\left( \xi ^{^{\prime }}(r)\right) ^{k}%
\right] ^{\prime }=r^{N-1}g^{k}(\xi \left( r\right) )\text{ in }(0,\rho ),
\\ 
\alpha _{n}=\xi _{n}\left( 0\right) =u_{n}\left( 0\right) \text{, }\xi
_{n}^{\prime }\left( 0\right) =0.%
\end{array}%
\right.  \label{mae}
\end{equation}%
then the solution $\xi _{n}\left( r\right) :=u_{n}\left( x\right) $ for $%
r=\left\vert x\right\vert $ of (\ref{mae}) can be extended so that remains a
solution of (\ref{mae}) in $\left[ 0,\rho _{0}\right) $. Next, apply Lemma %
\ref{2.4C4} with $\rho _{1}=\rho _{0}/2$ and $\rho _{2}=\rho _{0}$ to get%
\begin{equation}
\mathcal{K}\left( \beta _{n}\right) \geq \int_{\xi (\rho _{1})}^{\xi (\rho
_{2})}\frac{\left( C_{N-1}^{k-1}\right) ^{1/\left( k+1\right) }}{\sqrt[k+1]{%
\left( k+1\right) G\left( r\right) -G\left( \xi (\rho _{1})\right) }}dr\geq 
\frac{k}{N-2k}\rho _{0}^{\frac{2k}{k+1}}\left[ 1-\left( \frac{1}{2}\right) ^{%
\frac{N}{k}-2}\right] ,  \label{nd2k}
\end{equation}%
when\textit{\ }$N\neq 2k$\textit{\ and }%
\begin{equation}
\mathcal{K}\left( \beta _{n}\right) \geq \int_{\xi (\rho _{1})}^{\xi (\rho
_{2})}\frac{\left( C_{N-1}^{k-1}\right) ^{1/\left( k+1\right) }}{\sqrt[k+1]{%
\left( k+1\right) G\left( r\right) -G\left( \xi (\rho _{1})\right) }}dr\geq
\left( \frac{\rho }{2}\right) ^{\frac{2k}{k+1}}\ln 2,  \label{ne2k}
\end{equation}%
\textit{if }$N=2k$\textit{. }Passing to the limit as $n\rightarrow \infty $
in (\ref{nd2k})-(\ref{ne2k}), we obtain 
\begin{equation*}
0=\underset{\beta _{n}\rightarrow \infty }{\lim }\text{\ }\inf \mathcal{K}%
\left( \beta _{n}\right) \geq \left\{ 
\begin{array}{lll}
\frac{k}{N-2k}\rho _{0}^{\frac{2k}{k+1}}\left[ 1-\left( \frac{1}{2}\right) ^{%
\frac{N}{k}-2}\right] & \text{if} & N\neq 2k, \\ 
\left( \frac{\rho }{2}\right) ^{\frac{2k}{k+1}}\ln 2 & \text{if} & N=2k,%
\end{array}%
\right.
\end{equation*}%
which is a contradiction.

After these preliminaries we can begin to analyze the problem (\ref{ma1}).

\section{The Proof of Theorem \protect\ref{1.C4}\label{33}}

We are devoting this section to prove Theorem \ref{1.C4}. A natural way to
construct solutions in the case $k=1$ is by solving the finite datum
Dirichlet problem and then letting the datum grow to infinity to obtain the
conclusion (see also \cite{H}).

Indeed, let $u_{n}$ be the unique solution of the problem%
\begin{equation}
\left\{ 
\begin{array}{l}
\sigma _{k}^{1/k}\left( \lambda \left( D^{2}u\right) \right) =g(u)\text{ 
\textit{in}}\;\Omega , \\ 
\text{ }u_{\left\vert \partial \Omega \right. }=n,\text{ }\ n\in \mathbb{N}%
^{\ast },%
\end{array}%
\right.  \label{prob}
\end{equation}%
which clearly exists from Lemma \ref{min} in connection with Lemma \ref%
{subsuper} ii), for $\underline{u}=c=n$. Then $\underline{u}\leq u_{n}\leq 
\overline{u}$, where we remeamber that $\underline{u}$ is the solution of 
\begin{equation*}
\left\{ 
\begin{array}{l}
\sigma _{k}^{1/k}\left( \lambda \left( D^{2}u\right) \right) =g(n)+1\text{ 
\textit{in}}\;\Omega , \\ 
\text{ }u_{\left\vert \partial \Omega \right. }=n,\text{ }\ n\in \mathbb{N}%
^{\ast },%
\end{array}%
\right.
\end{equation*}%
and $\overline{u}$ is the solution of%
\begin{equation*}
\left\{ 
\begin{array}{l}
\Delta u=N\left( C_{N}^{k}\right) ^{-1}g\left( u\left( x\right) \right) 
\text{, }x\in B_{\tau }(x), \\ 
\underset{x\rightarrow x_{0}}{\lim }u\left( x\right) =+\infty \text{, }%
\forall \text{ }x_{0}\in \partial B_{\tau }(x).%
\end{array}%
\right.
\end{equation*}%
Since $u_{n}\leq u_{n+1}$ on $\partial \Omega $ it follows from Lemma \ref%
{comp} that $\{u_{n}\}_{n}$ is non-decreasing. Now, using the fact that $%
\partial \Omega $ is $C^{4+\alpha }$ it follows that there exists a radius $%
\tau >0$ such that for every $x_{0}\in \partial \Omega $ we can find a ball
of radius $\tau $, contained in $\Omega $ such that $x_{0}$ is also a
boundary point of the ball. Pick any boundary point $x_{0}\in \partial
\Omega $ and let $B_{\tau }(x)\subset \Omega $ be an interior ball
associated with $x_{0}$. The Lemma \ref{bila} shows that 
\begin{equation*}
\left\{ 
\begin{array}{lll}
\sigma _{k}^{1/k}\left( \lambda \left( D^{2}u\right) \right) =g\left(
u\right) & \text{in }B_{\tau }(x), &  \\ 
\underset{x\rightarrow x_{0}}{\lim }u\left( x\right) =+\infty & \forall 
\text{ }x_{0}\in \partial B_{\tau }(x), & 
\end{array}%
\right.
\end{equation*}%
has at least one blow-up solution $u_{\tau }$. By construction the solution
has radial symmetry. As a consequence of Lemma \ref{comp}, we can completely
answer the existence question for solutions of (\ref{ma1}). Indeed, we have
that 
\begin{equation*}
\underline{u}\left( x\right) \leq u_{n}\left( x\right) \leq u_{\tau }\left(
x\right) \text{ for }x\in B_{\tau }(x).
\end{equation*}%
This entails an upper bound for $u_{n}$. In particular, the sequence $u_{n}$
is uniformly bounded from above in $B_{\tau /2}(x)$. Notice that any $u_{n}$
is $k-$admissible. We next show that the sequence $\left\{ u_{n}\right\}
_{n} $ is uniformly bounded from above on every compact set included in $%
\Omega $. Pick any compact subset $\mathbb{K}\Subset \Omega $. Covering $%
\mathbb{K}$ by finitely many balls $B\left( x_{i},r_{i}/2\right) $ we
conclude that the sequence $\left\{ u_{n}\right\} _{n}$ is uniformly bounded
in $\mathbb{K}$, which is a compact set. Finally, since $\mathbb{K}$ is
arbitrary chosen it is clear that the limit%
\begin{equation*}
\lim_{n\rightarrow \infty }u_{n}\left( x\right) =u\left( x\right)
\end{equation*}%
exists as a continuous function and is a solution of 
\begin{equation*}
\begin{array}{ll}
\sigma _{k}^{1/k}\left( \lambda \left( D^{2}u\right) \right) =g\left(
u\right) & \text{in }\Omega ,%
\end{array}%
\end{equation*}%
and therefore $u_{n}\left( x\right) \overset{n\rightarrow \infty }{%
\rightarrow }u\left( x\right) $ on any compact subset $\mathbb{K}\Subset
\Omega $. There remains to prove that $u$ blows up at the boundary. This
question is clearly explained, for \ example in \cite{COV}. Let us sketch
the proof. We consider $x_{0}\in \partial \Omega $ and any sequence $x_{k}$
in $\Omega $ with $x_{k}$ $\rightarrow x_{0}$ as $k\rightarrow \infty $.
Since $\lim_{k\rightarrow \infty }u_{n}\left( x_{k}\right) =u_{n}\left(
x_{0}\right) =n$ on $\partial \Omega $ there is some $K>0$ such that $%
u\left( x_{k}\right) \geq u_{n}\left( x_{k}\right) \geq n$ when $k\geq K$.
Hence $\lim \inf_{k\rightarrow \infty }u\left( x_{k}\right) \geq \lim
\inf_{k\rightarrow \infty }u_{n}\left( x_{k}\right) =n$ and letting $n$ go
to infinity we can get that $u$ is a boundary blow-up solution of (\ref{ma1}%
) in $\Omega $.

It only remains to prove the reverse implication. For $n\in \mathbb{N}^{\ast
}$, we assume that $u_{n}$ solve (\ref{ma1}) in a ball $B$ of centre zero
and radius $\varepsilon _{n}$, where $\varepsilon _{n}$ is a decreasing
sequence such that $\varepsilon _{n}\searrow 0$ as\ $n\rightarrow \infty $.
Proceeding as in \cite{BAOII}, we can assume that $u_{n}$ is a radial
symmetric solution. Let $\beta _{n}=u_{n}(0)$. We observe that we can assume 
$\underset{n\rightarrow \infty }{\lim }\beta _{n}=\infty $, eventually by
passing to a subsequence. Then (\ref{ma9}) follows from (\ref{ma17}) applied
with $\rho =\varepsilon _{n}$. We finally prove that the sequence $\left\{
\beta _{n}\right\} $ is unbounded. If not, up to a subsequence, ($\beta _{n}$%
) converges to some $\beta \geq 0$. By (\ref{ma17}) applied with $\rho
=\varepsilon _{n}$, we have $u_{n}(\varepsilon _{n})=\infty $ and 
\begin{equation}
\begin{tabular}{l}
$0\leq \int_{u_{n}(0)}^{\infty }\frac{1}{\sqrt[k+1]{\left( k+1\right) \left(
G(z)-G(\beta _{n})\right) }}dz=\int_{0}^{\varepsilon _{n}}\frac{%
u_{n}^{^{\prime }}(r)}{\sqrt[k+1]{\left( k+1\right) \left(
G(u_{n}(r))-G(u_{n}(0))\right) }}dr$ \\ 
$\leq \frac{\left( C_{N-1}^{k-1}\right) ^{-1/\left( k+1\right) }\left(
k+1\right) }{2k}\varepsilon _{n}^{\frac{2k}{k+1}}.$%
\end{tabular}
\label{ult}
\end{equation}%
Passing to the limit as $n\rightarrow \infty $\ in (\ref{ult}) leads to%
\begin{equation*}
\int_{\beta }^{\infty }\frac{1}{\sqrt[k+1]{G(z)-G(\beta )}}dz=0
\end{equation*}%
which is not possible, since $\mathcal{K}(\beta )\in \left( 0,\infty \right) 
$.

\begin{remark}
Assume that $\psi $ belongs to a wide class $\Psi $ of monotone increasing
convex functions. There is an area in probability theory where
boundary-blow-up problems 
\begin{equation*}
\left\{ 
\begin{array}{l}
\Delta u=\psi \left( u\right) \text{ in }\Omega \text{ } \\ 
u=\infty \text{ on }\partial \Omega%
\end{array}%
\right.
\end{equation*}%
arise (see the paper \cite{din1} or directly the book \cite{din2} for
details). The area is known as the theory of superdiffusions, a theory which
provides a mathematical model of a random evolution of a cloud of particles.
Indeed, given any bounded open set $\Omega $ in the $N$-dimensional
Euclidean space, and any finite measure $\mu $ we may associate with these
the exit measure from $\Omega $ i.e. $\left( X_{\Omega },P_{\mu }\right) $,
a random measure which can be constructed by a passage to the limit from a
particles system. Particles perform independently $\Delta $-diffusions and
they produce, at their death time, a random offspring (cf. \cite{din3}). $%
P_{\mu }$ is a probability measure determined by the initial mass
distribution $\mu $ of the offspring and $X_{\Omega }$ corresponds to the
instantaneous mass distribution of the random evolution cloud. Then
proceeding in this way, one can obtain any function $\psi $ from a subclass $%
\Psi _{0}$ of $\Psi $ which contains $u^{\gamma }$ with $1<\gamma \leq 2$.
Dynkin \cite{din1}, also provided a simple probabilistic representation of
the solution for the class of problems $u^{\gamma }$ ($1<\gamma \leq 2$), in
terms of the so-called exit measure of the associated superprocess.
Moreover, the author says that a probabilistic interpretation is known only
for $1<\gamma \leq 2$.
\end{remark}

\begin{remark}
The problem of complex Hessian can be easily attacked (see \cite{X} as a
starting reference).
\end{remark}

\textbf{Acknowledgement. }I would like to thank to the editors and reviewers
for their advice and guidance in the review process.

\end{document}